\documentclass[reqno,a4paper]{amsart}
\usepackage[backend=biber,
style=numeric]{biblatex}
\usepackage{amsmath,amssymb,amsfonts,mathtools}
\usepackage[english]{babel}
\usepackage{bbold}
\usepackage{bigints}
\usepackage{pifont}
\usepackage[makeroom]{cancel}

\newtheorem{theorem}{Theorem}[section]

\newtheorem{corollary}[theorem]{Corollary}
\theoremstyle{definition}

\theoremstyle{definition}
\newtheorem{remark}[theorem]{Remark}
\theoremstyle{remark}
\newtheorem{example}[theorem]{Example}
\theoremstyle{remark}

\theoremstyle{remark}

\newcommand{\beq}{\begin{equation}}
\newcommand{\eeq}{\end{equation}}
\newcommand{\bdima}{\begin{displaymath}}
\newcommand{\edima}{\end{displaymath}}

\def\R{\mathbb R}

\newcommand{\I}{\textrm{I}}

\begin{document}

\title[On the Ultimate Boundedness]{On Spatial Cohesiveness of Second-Order Self-Propelled Swarming Systems}

\author[C. Medynets]{Constantine Medynets}
\email{medynets@usna.edu}
\address{Mathematics Department, United States Naval Academy, Annapolis, MD 21402}

\author[I. Popovici]{Irina Popovici}
\email{popovici@usna.edu (corresponding author)}
\address{Mathematics Department, United States Naval Academy, Annapolis, MD 21402}

\thanks{The research of C.M. and I.P. was supported by Office of Naval Research Grant \# N0001421WX00045}


\begin{abstract} The study of emergent behavior of  swarms is of great interest for applied sciences. One of the most fundamental questions for self-organizing swarms is whether the swarms disperse or remain in a spatially cohesive configuration. In the paper we study  dissipativity properties and spatial cohesiveness of the swarm of self-propelled  particles governed by the model $\ddot r_k = -p_k(|\dot r_k|)\dot r_k - \sum_m a_{k,m}r_m$, where $r_k\in \mathbb R^d$, $k=1,\ldots,n$, and $A = \{a_{k,m}\}$ is a symmetric positive-semidefinie matrix. The self-propulsion term is assumed to be continuously differentiable and to grow faster than $1/z$, that is, $p_k(z)z\to\infty $ as $z\to\infty$. We establish that the velocity and acceleration of the particles are ultimately bounded. We show that when $\ker (A)$ is trivial, the positions of the particles are also ultimately bounded. For systems with  $\ker (A)\neq \{0\}$, we show that, while the system might infinitely drift away from its initial location, the particles remain within a bounded distance from the generalized center of mass of the system, which geometrically coincides with the weighted average of agent positions. The weights are determined by the coefficients of the projection matrix  onto $\ker (A)$.

In our proof we switch to the velocity-acceleration coordinates and focus on the study of dissipativity properties for a more general class of Li\'enard systems $\ddot x_k = -\mathbb F_k(x_k)\cdot \dot x_k -\sum_{m} a_{k,m}x_m$, $k=1,\ldots,n$, $\mathbb F_k(x) = \nabla F_k(x)$ with $F_k: \mathbb R^d\rightarrow \mathbb R^d$ given by $F_k(x) = p_k(|x|)x$. The phase space of the system splits into a family of invariant manifolds determined by the kernel of the matrix $A$. We establish that this system is ultimately bounded within each of these invariant manifolds. We also include the proof of the ultimate boundedness of velocities and accelerations for systems with bounded coupling, including systems coupled via the Morse potential. 
\end{abstract}

\maketitle

\section{Introduction}

The study of emerging collective behavior in biological or man-made systems has been an active area of interdisciplinary research drawing the interest of engineers, applied scientists, mathematicians, physicists, and many others.    One of the most prominent and defining  features of any biological swarm is its ability to form spatially coherent configurations. Using the language of dynamical systems, we can say that the models that describe biological systems must be inherently relatively dissipative  and lead to the congregation of the swarm around its  center. Due to large dimensionality of any biological system, their rigorous mathematical studies involving detailed proofs are rare.  The present paper deals with a mathematical study of a large  class of {\it second-order} swarming models that appear in various applications to biological and robotic systems  \cite{Mach:2006,Szwaykowska:2016,Ebeling:2001,EdwardsSchwartz:2020}.

The present paper has grown out of our interest in the study of dynamical properties of the swarming model
\begin{equation}\ddot r_k = (1-|\dot r_k|^2)\dot r_k-\frac{1}{n}\sum_{m=1}^n(r_k- r_m),
\label{theequation}
\end{equation} where $r_k\in \mathbb R^2$ represents the position vector of the $k$-th agent in the swarm and $\dot r_k$ stands for the time derivative $\frac{dr_k}{dt}$. Heuristically, the model presented by Equation (\ref{theequation}) describes the system in which each agent tries to maintain a unit speed in the direction it has been traveling and at the same time there is a spring-like attraction force between each pair of agents.

In  \cite{Szwaykowska:2016,EdwardsSchwartz:2020} the authors have conducted a series of mixed-reality experiments involving robotic cars and boats whose motion and interactions are described by Equation (\ref{theequation}). We should point out that the authors consider the Langevin version of the model with added noise and a repulsion term.   Equation (\ref{theequation}) also appears in the study of biological systems \cite{Mach:2006},  where it was used to model vortex swarming in crustaceans. More on modeling of circle swimmers as well as a literature review can be found in \cite{Kaiser:2013}.  In literature, Model (\ref{theequation}) is often referred to as the {\it parabolic potential model}.

Most differential-equation-based second-order swarming models follow the pattern: 
$$\ddot r_k = \textrm{self-propulsion} + \textrm{coupling},$$
The self-propulsion term describes the way each agent gains kinetic energy from their environment or from within. The coupling term describes mutual attraction and repulsion between agents. Most {\it numerical experiments} on swarming models, see, for example, \cite{Romanczuk:Book,DorsognaChuangBertozziChayes:2006,Ebeling:2001}, show that  the velocity and acceleration of any solution $r(t)$, regardless of the initial conditions, are {\it ultimately bounded}, in the sense that there is a universal constant $C>0$ that only depends on the number of agents $n$ such that  $|\dot r_k(t)|\leq C$ and $|\ddot r_k(t)|\leq C$ for all large enough $t$.   The ultimate boundedness of the velocity and acceleration in Model (\ref{theequation}) automatically implies  that there is a constant $K>0$ that only depends on $n$ such that $|r_k-\bar r|\leq K$ as $t\to\infty$, where $\bar r = (r_1+\cdots+r_n)/n$. In nature swarms often converge to  spatially tight configurations around their center of mass, which is also a desirable feature for many physical systems. 

The main goal of the present paper is to  establish the ultimate boundedness (dissapitivity) of the velocity and acceleration (and when applicable, the position) for the class of swarming models of self-propelled particles  moving in $\R^d$ whose position vectors $r_k(t) \in \R^d $  satisfy

\begin{equation}\label{EqnGeneralModel}
\ddot r_k = -p_k(|\dot r_k|) \dot r_k -\sum_{m=1}^n a_{k,m}r_m, \;  k=1,\ldots,n
, \end{equation}
where $A = \{a_{i,j}\}$ is a  positive semidefinite symmetric matrix. We also prove the ultimate boundedness for dynamical systems governed by the equation  

\begin{equation}\label{EqnGeneralModelBounded}
\ddot r_k = -p_k(|\dot r_k|) \dot r_k +c_k(r_1,\ldots,r_n,\dot r_1,\ldots, \dot r_n), \;  k=1,\ldots,n
, \end{equation}
under the assumption that the coupling functions $c_k$, $k=1,\ldots,n$, which might depend on both the position and the velocity of the agents in the swarm, are uniformly bounded.  In both cases, the functions  $p_k(z)$, $k=1,\ldots,n$, are assumed to be continuously differentiable and satisfy $p_k(z)z\to\infty$ as $z\to\infty$. Note that Morse potential  systems  \cite{DorsognaChuangBertozziChayes:2006,Armbruster:2017} fall under the category of systems with bounded coupling.

We observe that Equation (\ref{theequation}) is obtained from  (\ref{EqnGeneralModel}) by setting $A = \mathbb I - Q$, where $\mathbb I$ is the identity matrix and  
 $$Q = \frac{1}{n}\left(\begin{array}{cccc} 1 & 1 & \cdots & 1 \\
 1 & 1 & \cdots & 1 \\
 \vdots & \vdots & \ddots & \vdots \\
 1 & 1 & \cdots & 1  \end{array}\right).$$

 In recent years there has been a significant interest in understanding the collective behavior of swarms consisting of non-identical agents, and of swarms whose communication topology is such that each member only senses the position of a fixed number of initial neighbors. In particular, recent studies include observing birds in natural settings (flocks of hundred of Surf Scoters  in
\cite{Lukeman12576}
 or thousands of European Starlings in \cite{Ballerini1232}), where it was empirically established that each bird's motion is controlled by at most seven neighbors,  biomimetric explorations (autonomous mobile robots or air vehicles whose proximity sensors communicate with only two nearest neighbors in
\cite{1166526}), and hydrodynamic models describing the
collective behavior of a continuum of agents
\cite{Shu_2020}. Switching from all-to-all coupling in (\ref{theequation}) to weighted coupling in (\ref{EqnGeneralModel}) allows us to model different communication topologies in the swarm. Choosing different propulsion terms in (\ref{EqnGeneralModel}) also allows to consider heterogeneous swarms. 

\begin{figure}[h!]
\centering
\includegraphics[width=.3\textwidth]{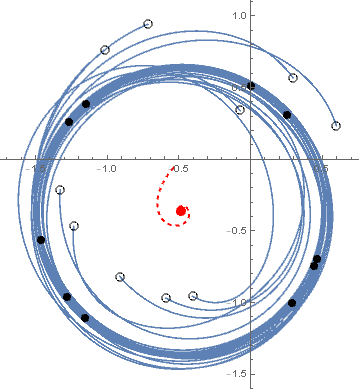}
\caption{\label{FigureRingState} Credit \cite{kolon2021stability}. Asymptotic behavior of (\ref{theequation}) in the basin of attraction of a non-degenerate ring state of the parabolic potential model. The center of mass $\bar r(t)$  is indicated as a red dot. The dashed line represents the trajectory of $\bar r(t)$.  The  initial (terminal) positions of the agents at  are represented as hollow (solid) dots.}
\end{figure}

While running simulations on Equation (\ref{EqnGeneralModel}), one can notice that swarms often congregate around certain points in space. 
For example, for the parabolic potential model, the center of mass is such a point, leading to the onset of  the {\it ring/milling state} --  when the center of mass of the swarm is stationary and all agents rotate around it with a unit speed, see Figure (\ref{FigureRingState}). Thus, the geometric center of the swarm plays a significant role in the pattern formation of systems governed by Equation (\ref{theequation}).  We would like to mention that the study of the stability of ring states has recently received a lot of attention, see, for example, \cite{Szwaykowska:2016} and references therein for numerical studies and \cite{kolon2021stability} for detailed proofs. 
  
For a general swarm system satisfying Equation (\ref{EqnGeneralModel}), in Corollary \ref{CorollaryParabolicPotential3} we will show that  the coordinate-wise projection of the corresponding coordinates of $\{r_1,\ldots,r_n\}$  onto the kernel of the matrix $A$ plays the role of the generalized center of mass of the swarm, whereby finding the precise mechanism for the occurrence of special ``gravitational'' centers.  We note that the center of mass of the swarm governed by the parabolic potential model is precisely the coordinate-wise projection of the corresponding coordinates of the agent positions onto the kernel of the coupling matrix $A$. A similar situation also holds for many systems with sparse communication topologies, say, with symmetric nearest neighbor coupling, see Example \ref{ExamplesGeomCenters} for details.

Although our results are formulated in the language of multi-agent ODEs and swarms, they have broader reach. For example, the functions $r_k(t)$ can represent the spatial discretization of the solution to a damped nonlinear wave equation, with the positive definite matrix $A$ being the finite difference approximation  of the Laplacian \cite{Hale:1997}. Establishing the existence of uniform bounds (which holds in the case of  coupling with a nonsingular matrix $A$) for the position and velocity vectors is a required first step in understanding the global attractor for the system.

Switching to the velocity-acceleration coordinates (setting $x_k=\dot r_k$) in System (\ref{EqnGeneralModel}) with linear coupling, we can rewrite (\ref{EqnGeneralModel}) as the system of coupled Li\'enard oscillators
\begin{equation}\label{eqnLienardSystemIntro}\ddot x_k = -\mathbb F_k(x_k)\cdot \dot x_k -\sum_m a_{k,m} x_m, \; k=1,\ldots,n\end{equation}
 where  $\mathbb F_k(x) = \nabla F_k(x)$ and $F_k(x) = p_k(|x|)x$, see the details in Section \ref{SectionLinearCoupling}. We notice that in the case when $d=1$, ($x_k\in \mathbb R$), in (\ref{eqnLienardSystemIntro}),  we obtain a system of coupled van der Pol equations.

 We would like to mention several mathematical papers in which systems similar to Equation (\ref{theequation}) have been studied. General second-order gradient-like systems $\ddot r = f(\dot r)-  \nabla U(r)$ were studied in \cite{Haraux:1998} and \cite[Ch. 7]{Haraux:book}, where it was shown that under some conditions on $f$ and the potential function $U$, most notably $f(\dot r)\cdot \dot r<0$ if $\dot r\neq 0$, every bounded solution converges to a configuration that solves $\nabla U(r)=0$. We note that the functions $-p_k(|\dot r_k|)||\dot r_k|^2$ in (\ref{EqnGeneralModel}) need not be negative-definite and, thus, the results of \cite{Haraux:1998} do not apply.  Not imposing that $p_k>0 $ places our results outside of the purview of \cite{Tunc:2018} (specifically their A2 assumption).

 A comprehensive discussion of the asymptotic behavior
of coupled dissipative systems whose coupling strength is amplified to infinity can be found in \cite{Hale:1997}. Moreover, owing to several assumptions placed on the nonlinearities in first order systems (2.17 in  \cite{Hale:1997}), and on nonlinearities for second order systems (4.1 in \cite{Hale:1997}), our results are not encompassed by \cite{Hale:1997}.

 Our approach has been informed by the ideas from \cite{Figueiredo:Boundedness}, where the authors considered  a boundedness condition for coupled 1D systems, see also a related discussion in  \cite{Skoog:1972}.  We note that the proof in  \cite{Figueiredo:Boundedness} is incomplete and it is unclear whether their proof can be rectified based on the assumptions made in the paper. The ultimate boundedness of some coupled 1D systems is also discussed in \cite{Kato:1986}.  We note that the assumptions made in \cite{Figueiredo:Boundedness} or \cite{Kato:1986} do not apply to the systems discussed in the present paper.

In Theorem \ref{ThmMain}  we establish the ultimate boundedness of the velocity-acceleration coordinates for systems  (\ref{EqnGeneralModel}) with linear coupling. In Corollary \ref{CorollaryParabolicPotential2} we show that systems  (\ref{EqnGeneralModel}) with linear coupling are ultimately bounded whenever the coupling matrix is invertible. In Corollary \ref{CorollaryParabolicPotential} we obtain the ultimate boundedness result for the parabolic potential model (\ref{theequation}). Theorem \ref{TheoremWeakCoupling} establishes the ultimate boundedness for the velocities and accelerations of swarms with bounded coupling.
We note that unlike systems with linear coupling  systems with bounded coupling may have positions that diverge away from the center of mass  \cite{DorsognaChuangBertozziChayes:2006}.

{\bf Acknowledgement.} We learned about System (\ref{theequation}) and  applications of swarming systems to robotics from Dr. Ira Schwartz, Dr. Jason Hindes, and Dr. Klementyna Szwaykowska of the Naval Research Laboratory, Washington, D.C., when they presented their research in the U.S. Naval Academy Applied Mathematics Seminar. C. Medynets and C. Kolon also visited Dr. Schwartz and his research group at the NRL in the summer of 2016 and 2017, respectively. We are thankful for their hospitality and fruitful discussions. C. M. and I.P. also acknowledge the support from the Office of Naval Research, Grant \# N0001421WX00045.

%
\section{Linear Coupling}\label{SectionLinearCoupling}

 For a (column) vector-valued function $F=(F_1,\ldots, F_d)^T$, $F_i: \mathbb R^d\rightarrow \mathbb R$, we will write $\nabla F(x)$ to denote the matrix
$$\displaystyle \nabla F(x) = \left (\begin{array}{cccc}
\vspace{0.05in}
\frac{\partial F_1}{\partial x_1} & \frac{\partial F_1}{\partial x_2} & \cdots & \frac{\partial F_1}{\partial x_d} \\
\frac{\partial F_2}{\partial x_1} & \frac{\partial F_2}{\partial x_2} & \cdots & \frac{\partial F_2}{\partial x_d} \\
\vdots & \vdots & \ddots & \vdots \\
\frac{\partial F_d}{\partial x_1} & \frac{\partial F_d}{\partial x_2} & \cdots & \frac{\partial F_d}{\partial x_d}
 \end{array}\right).$$  In what follows, all vectors are assumed to be column-vectors. By $'\cdot'$ we will denote the regular matrix multiplication. For a matrix/vector $v$, $v^T$ will stand for its transpose. We will denote the standard Euclidian norm of  $v\in \mathbb R^d$ by $|v|$. For a matrix $L$, the $L^{\infty}$ operator norm, equal to the maximum absolute row sum of the matrix, will be denoted by
 $||L||_{\infty}$.

In the summations that follow, the summation index is assumed to run through the index set $\{1,\ldots,n\}$  unless stated otherwise.

 Consider the dynamical system given by

\begin{equation}\label{eqnMainSystem} \left\{
\begin{array}{l}
\dot x_k  = y_k  \\
\dot y_k  = - \mathbb F_k (x_k)\cdot y_k-\sum\limits_m a_{k,m} x_m
\end{array}
\right.
\end{equation}
where $x_k,y_k\in \mathbb R^d$, $k=1,\ldots,n$, the matrix $A = \{a_{i,j}\}$ is symmetric and positive semidefinite, $\mathbb F_k(x) = \nabla F_k(x)$,  $F_k(x) = p_k(|x|)x$ and $p_k : [0,\infty)\rightarrow \mathbb R$ is a continuously differential function that satisfies $p_k(z)z \to \infty$  as $z\to\infty$.

  Denote by $Q = \{q_{i,j}\}$ the matrix of the orthogonal projection onto the kernel of the matrix $A$. If the matrix $A$ is invertible, then $Q$ is the zero matrix. The projection $Q$ can be constructed as $Q = Q^T=K K^T $, where $ K$ is an $n\times \textrm{dim(Ker(A))}$-matrix whose columns form an orthonormal basis for the kernel of $A$. It follows from the fact that $A^T=A$ and $AQ=0$ that  $QA = 0$.  Therefore, $\sum_{l} q_{l,k}a_{k,m} = 0$ for any $l,m\in \{1,\ldots, n\}$. Thus, multiplying both sides of (\ref{eqnMainSystem}) by $Q$, we obtain that
  $$\sum\limits_k q_{l,k}\left [\dot y_k  + \mathbb F_k(x_k)\cdot y_k \right ]
  = \sum_{k}q_{l,k}\sum_m a_{k,m} x_m = \sum_{m}\left( \sum_k q_{l,k} a_{k,m}\right)  x_m  = 0_d,$$ the $d$-dimensional zero vector.
  Thus, for every $1\leq l \leq n$, we obtain that
\begin{equation}\label{eqnEnergy}
\sum_kq_{l,k}\left[ y_k + F_k(x_k) \right] = E_l =  \textrm{constant}
\end{equation}
and, therefore, this vector equation defines $d$ invariant manifolds in the phase space. We will refer to the function defined in (\ref{eqnEnergy}) as the {\it energy} of the system. Note that some rows of $(\ref {eqnEnergy} )$  are linear combinations of the others, resulting in fewer than $nd$ algebraically independent invariant manifolds.  In the following result we will show that the system is dissipative within each constant energy manifold defined by (\ref{eqnEnergy}).

 %

\begin{theorem}\label{ThmMain}  For any collection of vectors $E_1,\ldots, E_n\in \mathbb R^d$ for which the invariant manifold $\mathcal M$ determined by $\sum_kq_{l,k}\left[ y_k + F_k(x_k) \right] = E_l$, $1\leq l\leq n$, is non-empty, there are constants $a>0$ and $b>0$  such that for any solution $(x_k,y_k)$ of (\ref{eqnMainSystem}) with initial conditions lying in the invariant manifold $\mathcal M$, we have that $|x_k(t)|\leq a$ and $|y_k(t)|\leq b$ for all $t$ large enough.
\end{theorem}
{\it Proof.} (I)  Introduce the auxiliary ``ramp'' function $M:\mathbb R^d \rightarrow \mathbb R^d$ as follows:
$$M(x) =
\left\{\begin{array}{cc}
x & \mbox{if } |x|\leq 1 \\
\frac{x}{|x|} & \mbox{if } |x|>1
 \end{array}\right.
 $$
The function $M(x)$ is globally Lipschitz continuous and it is smooth everywhere but $|x| =1$.  Using a slight abuse of notation, define
\begin{equation}\label{eqnGradientM}
\nabla M(x) =
\left\{\begin{array}{cc}
\I_d & \mbox{if }|x|\leq 1 \\
\frac{1}{|x|^3}\left[|x|^2\I_d - (x\cdot x^T)  \right]  & \mbox{if } |x|>1
\end{array}\right.
\end{equation}
 where $\I_d$ is the $d\times d$ identity matrix. Notice that $\nabla M(x)$ coincides with the (conventional) gradient matrix of $M$ everywhere but $|x| = 1$ where the gradient does not exist in the classical sense.

 Consider the matrix $S = A+Q$ and denote its entries by $S = \{s_{i,j}\}$.
 Since the matrix $A$ is symmetric, in view of the spectral theorem,  it can be represented as $A = \sum_{\lambda \in \textrm{Sp(A)}}\lambda Q_\lambda$, where $Q_\lambda$ is an orthogonal projection onto the eigenspace corresponding to the eigenvalue $\lambda$. Furthermore, these projections and eigenspaces are mutually orthogonal, which implies that $S$ has the same eigenspaces and the corresponding eigenvalues as the matrix $A$ except for the vectors in the null-space of $A$ that become eigenvectors of  $S$ corresponding to the eigenvalue $\lambda =1$. In particular,  $S$ is symmetric, positive-definite, and invertible.


 Denote by $\{s_{i,j}^{(-1)}\}$ the coefficients of the matrix $S^{-1}$ and by $\{\epsilon_{i,j}\}$ the coefficients of the identity matrix $I$. Note  that $S(I-Q) = (A+Q)(I-Q)=A-AQ+Q-Q^2=A $ since $AQ=0$ and $Q^2=Q.$
 Therefore,

 \begin{equation}\label{eqnMatrixIdentity} S^{-1}A = I-Q.  \end{equation}

Consider the function
$$V = \frac{1}{2}\sum_k x_k^T\cdot x_k + \frac{1}{2}\sum_{k,m}s_{k,m}^{(-1)} \left[y_k+ F_k(x_k)-\delta W_k\right]^T\cdot \left[y_m+ F_m(x_m)-\delta W_m\right],$$
where $a$ are $\delta$ are some real numbers, with $a>1$,
 $0<\delta <1 $, and
$$W_k =\sum_{l}  s_{k,l} M\left(\frac{x_l}{a}\right).$$
Since the matrix $S^{-1}$ is positive-definite, the function $V$ is non-negative with $V>0$ for all points except for the origin in $\mathbb R^{2nd}.$

 Our goal is to show that for some choice of  the parameters $a>0$ (large), $0<\delta<1$ (small),  and $b>0$ (large), the function $V$ decreases at a rate $\dot V\leq -1 $ along the trajectories of (\ref{eqnMainSystem}) in the complement of the box $\{(x_k,y_k) : |x_k|\leq a,\; |y_k|\leq b,\;k=1,\ldots,n\}$.  The parameters $a$ and $\delta$ will be chosen so that  that $\dot V\leq-1$ whenever $|x_k|>a$ for at least one $k$. The choice of the parameter $b$ will guarantee that $\dot V\leq -1$ whenever $|x_k|\leq a$ for all $k$ and $|y_k|>b$ for at least one $k$. The result will then follow from an application of the Lyapunov method
 \cite{LaSalle:Book:1961}.

We note that the function $V$ is continuous and it is also smooth everywhere but the points with $|x_k| = a$. Thus, we have to be extra careful in our treatment of the derivative of  $V$ along the trajectories of the system and we need to carefully explain in what sense this function is differentiable.  Let $\gamma(t) = (x(t),y(t))$ be a $C^2$ solution of  (\ref{eqnMainSystem}), where $x = (x_1,\ldots,x_n)$ and $y = (y_1,\ldots,y_n)$.  Consider the function $t\mapsto V(\gamma(t))$. Note that this function is locally Lipschitz. Denote by $\dot V$ the {\it weak} derivative of this function. Our first objective is to justify that the weak derivative $
\dot V$ can be calculated using the standard differentiation rules and substituting (\ref{eqnGradientM}) for the gradient of $M$.

Recall that for smooth functions, the weak derivative coincides with the conventional derivative. Expanding the function $V$, we can represent it as the linear combination of pairwise dot products involving the functions $x_k$, $F_k(x_k)$, $y_k$, and $M\left(\frac{x_k}{a} \right)$.
Note that the  terms that do not contain the function $M\left(\frac{x_k}{a} \right)$ are smooth functions and, thus, can be differentiated using the standard product and chain rules for multivariate functions.  The derivative of the terms involving  the function $M\left( \frac{x_k}{a}\right)$  has to be evaluated using the product rule for  weak derivatives,  followed by the standard chain rule for $y_k$ and $F_k(x_k)$ and the the chain rule for weak derivatives of $M\left( \frac{x_k}{a}\right)$. We note that one has to be careful when applying the chain rule for weak derivatives  as such a chain rule does not always hold. An interested reader can find a general discussion of the (weak derivative) chain rule for multivariate functions and some situations when it fails  in \cite[Section 4.3]{Leoni:Book}.  For our purposes, we  need to justify only why the weak derivative of $M(x(t))$ along the trajectory $\gamma(t)=(x(t),y(t))$ is equal to $\nabla M(x)\cdot y $.

We observe that the function $M(x)$ falls under the assumptions of Theorem 2.1 in \cite{MuratTrombetti:2003} since the function $M(x)$ is (1) globally Lipschitz-continuous,  (2) satisfies $M(0)=0$, and (3) constructed by two $C^1$ functions $M_1$ and $M_2$ given by $M_1: \{|x|\leq 1\}\rightarrow \mathbb R^d$, $M_1(x) =x$, and $M_2 : \{|x|>1\}\rightarrow \mathbb R^d$, $M_2(x)  =x/|x|$,  that admit Lipschitz-continuous and $C^1$ extensions to all of $\mathbb R^d$. Therefore, applying \cite[Theorem 2.1]{MuratTrombetti:2003} for any interval $[t_1,t_2]$, we obtain that the function $t \mapsto M(x(t))$ belongs  the Sobolev space $W^{1,\epsilon}([t_1,t_2],\mathbb R)$,  for all $\epsilon >1$, and its weak derivative equals
$$\frac{d M(x(t))}{dt} = \nabla M(x(t))\cdot y(t) \mbox{ for almost every }t.$$

Therefore, the function $t\mapsto V(\gamma(t))$ is differentiable almost everywhere and its weak derivative can be computed using the standard differentiation rules for almost all $t$. In what follows, when we talk about derivatives, we mean weak derivatives which exist for almost every $t$ and all the inequalities are assumed to hold for almost every $t$.

\medskip

(III) Note that the (full) derivative of $y_m+ F_m(x_m)-\delta W_m$ along the trajectories  of (\ref{eqnMainSystem})  is equal to $$-\mathbb F_m(x_m)\cdot y_m-\sum_l a_{m,l}x_l + \mathbb F_m(x_m)\cdot y_m- \delta \dot W_m = -\sum_l a_{m,l}x_l - \delta \dot W_m,$$ where $\dot W_m$ is the full derivative of $W_m$.

Using the fact that the matrix $S^{-1}$ is symmetric, we obtain that the full derivative of $V$ is
\begin{equation}\label{eqnDotVAux1}
\begin{split}
\dot V  &= \sum_k y_k^T\cdot x_k + \sum_{k,m}s_{k,m}^{(-1)}\left[y_k +F_k(x_k)-\delta W_k \right]^T
\left [ -\sum_l a_{m,l}x_l - \delta \dot W_m \right] \\
&=\!\begin{multlined}[t]
\sum_k y_k^T\cdot x_k \\ + \sum_{k}\left[y_k +F_k(x_k)-\delta W_k \right]^T \cdot
\left [ -\sum_l\left[\sum_m s_{k,m}^{(-1)} a_{m,l}\right] x_l - \delta \sum_m s_{k,m}^{(-1)} \dot W_m \right].
\end{multlined}
\end{split}
\end{equation}

 It follows from (\ref{eqnMatrixIdentity}) that
 $$\sum_l\sum_m s_{k,m}^{(-1)} a_{m,l}x_l  = \sum_l \left ( \epsilon_{k,l}x_l-q_{k,l}x_l \right ) = x_k - \sum_l q_{k,l}x_l.$$
 Note also that
 $$ \sum_m s_{k,m}^{(-1)} \dot W_m  = \sum_l\sum_{m} s_{k,m}^{(-1)}s_{m,l} \dot M\left(\frac{x_l}{a} \right) = \sum_l \epsilon_{k,l} \dot M\left(\frac{x_l}{a} \right) = \dot M\left(\frac{x_k}{a} \right).$$

 Substituting these equations into (\ref{eqnDotVAux1}), and distributing from right to left, we obtain that

 \begin{equation}\label{eqnDotVAux2}
\begin{split}
\dot V & =  \sum_k y_k^T\cdot x_k + \sum_{k}\left[y_k^T +F_k(x_k)^T-\delta W_k ^T \right] \cdot
\left [- x_k + \sum_l q_{k,l}x_l - \delta \dot M\left(\frac{x_k}{a} \right) \right] \\
& = \!\begin{multlined}[t]
-\delta \sum_k \left[y_k +F_k(x_k)-\delta W_k \right]^T\cdot \dot M\left(\frac{x_k}{a} \right)
 -\delta \sum_kW_k^T\cdot \left(-x_k+\sum_l q_{k,l}x_l \right)+ \\
+\sum_{l}\sum_{k}q_{k,l}\left[y_k +F_k(x_k)\right]^T\cdot x_l- \sum_k F_k(x_k)^T\cdot x_k
 \end{multlined}
\end{split}
\end{equation}

Recall that the matrix $Q$ is symmetric. It follows from (\ref{eqnEnergy}) that $$\sum_{l}\sum_{k}q_{k,l}\left[y_k +F_k(x_k)\right]^T\cdot x_l = \sum_l E_l^T\cdot x_l.$$ Thus, we can rewrite Equation (\ref{eqnDotVAux2}) as
 \begin{multline}\label{eqnDotVAux3}
\dot V  =    - \sum_k F_k(x_k)^T\cdot x_k +\sum_{k}E_k^T\cdot x_k -\delta \sum_kW_k^T\cdot \left(-x_k+\sum_l q_{k,l}x_l \right) \\
 -\delta \sum_k \left[y_k +F_k(x_k)-\delta W_k \right]^T\cdot \dot M\left(\frac{x_k}{a} \right)
\end{multline}

Fix an index $k_1$ such that $|x_{k_1}|= \max_k|x_k|$. Then, using the Cauchy-Schwartz inequality we obtain that
$$|\sum_k E_k^T\cdot x_k|\leq \left ( \sum_k |E_k| \right ) |x_{k_1}|.$$
Recall that the ramp function satisfies $ |M(x)|\leq 1 $ for all $x \in \mathbb R^d, $
so $|W_k|\leq \sum_l |s_{k,l}| \leq ||S||_{\infty} $.
It follows from the Cauchy-Schwartz inequality,  that
$$|\sum_k W_k^T\cdot (-x_k+\sum_l q_{k,l}x_l ) | \leq ||S||_{\infty} \sum_k ( |x_k|+ |\sum_l q_{k,l}x_l|) \leq n ||S||_{\infty}  \left(1+||Q||_{\infty} \right)|x_{k_1}|.$$

Combining the last two inequalities with Equation (\ref{eqnDotVAux3}) and using the fact that
$\dot M\left(\frac{x_k}{a} \right)=\frac{1}{a}\nabla M\left(\frac{x_k}{a}\right)\cdot y_k$,
we obtain that
\begin{displaymath}
\dot V  \leq  - F_{k_1}(x_{k_1})^T\cdot x_{k_1} + C_1|x_{k_1}| - \sum_{k\neq k_1} F_k(x_k)^T\cdot x_k+
\end{displaymath}
\begin{equation}\label{eqnDotVAux4}
 -\frac{\delta}{a} \sum_k \left[y_k +F_k(x_k)- \delta W_k \right]^T\cdot \nabla M\left(\frac{x_k}{a}\right)\cdot y_k,
\end{equation}
where
$C_1 = \left(\sum_k |E_k|\right)+  n ||S||_{\infty} (1+||Q||_{\infty})$.

We need to rearrange the last sum in the right-hand side of the previous inequality according to whether $|x_k|> a$. Using the definitions of $M(x/a)$, $\nabla M(x/a)$, and $F(x)$, we obtain that those terms sum up to
\begin{equation}\label{eqnZAux}
\begin{split} Z &=  - \frac{\delta}{a}\sum_{|x_k|>a} \left[y_k+ F_k(x_k)-\delta W_k\right]^T \cdot \nabla M\left(\frac{x_k}{a}\right)\cdot y_k  \\
 &=  - \frac{\delta}{a}\sum_{|x_k|>a} \left[\left[y_k+ F_k(x_k)\right]^T\cdot  \nabla M\left(\frac{x_k}{a}\right)\cdot y_k -\delta W_k^T \cdot \nabla M\left(\frac{x_k}{a}\right)\cdot y_k \right ]\\
& =   - \frac{\delta}{a}\sum_{|x_k|>a} \left[ a\left[y_k^T+ p_k(|x_k|)x_k^T\right]  \cdot \left(|x_k|^2\I_d-x_k\cdot x_k^T\right)\cdot \frac{y_k}{|x_k|^3}-\delta W_k^T\cdot \nabla M\left(\frac{x_k}{a}\right)\cdot y_k \right]\\
 & = \!\begin{multlined}[t] - \frac{\delta}{a}\sum_{|x_k|>a}\bigg[ \frac{a}{|x_k|^3}\Big [
 |x_k|^2 |y_k|^2-(y_k^T\cdot x_k)^2+p_k(|x_k|)|x_k|^2x_k^T\cdot y_k
 -p_k(|x_k|)|x_k|^2x_k^T\cdot y_k \Big] \\ -\delta W_k^T \cdot \nabla M\left(\frac{x_k}{a}\right)\cdot y_k \bigg]
  \end{multlined}
   \\
  & =   - \frac{\delta}{a}\sum_{|x_k|>a}\left[ \frac{a}{|x_k|^3}\Big [
 |x_k|^2 |y_k|^2-(y_k^T\cdot x_k)^2\Big] -\delta W_k^T \cdot \nabla M\left(\frac{x_k}{a}\right)\cdot y_k \right].
  \end{split}
\end{equation}
For vectors $x$ and $y$ in $\mathbb R^d$ denote by $\textrm{ort}_xy$ the vector rejection of $y$ from $x$, that is, $$\textrm{ort}_xy = y - \textrm{proj}_xy = y - \frac{y^T\cdot x}{|x|^2}x.$$
Notice that
$$|\textrm{ort}_xy|^2 = |y|^2-\frac{(y^T\cdot x)^2}{|x|^2} = \frac{|x|^2|y^2|-(y^T \cdot x)^2}{|x|^2}.$$
Notice also that
$$ \nabla M\left(\frac{x_k}{a}\right)\cdot y_k = \frac{a}{|x_k|}\left[y_k-\frac{x_k^T\cdot y_k}{|x_k|^2}x_k \right] = \frac{a}{|x_k|}\textrm{ort}_{x_k}y_k.$$

Thus, we can rewrite (\ref{eqnZAux}) as
\begin{equation}\label{eqnZAux2}
\begin{split}Z &= -\frac{\delta}{a}\sum_{|x_k|>a} \left [  \frac{a}{|x_k|}|\textrm{ort}_{x_k}y_k|^2-\frac{\delta a}{|x_k|} W_k^T\cdot \textrm{ort}_{x_k}y_k  \right ] \\
&= -\delta\sum_{|x_k|>a}  \frac{1}{|x_k|}\Big[ |\textrm{ort}_{x_k}y_k|^2-\delta \cos(\theta_k) |W_k| |\textrm{ort}_{x_k}y_k| \Big],
\end{split}
\end{equation}
where $\theta_k$ is the angle between $W_k$ and $\textrm{ort}_{x_k}y_k$. Notice that each term in (\ref{eqnZAux2}) is a quadratic function in $|\textrm{ort}_{x_k}y_k|$. Recall that any parabola of the form $z^2+ d z$ attains its global minimum at $z = -d/2$. Thus, $z^2+ d z\geq \frac{-d^2}{4}$ for any $z\in \mathbb R$. Therefore, since $|W_k|\leq ||S||_{\infty} $ and $|\cos(\theta_k)|\leq 1$, there is a constant $C_2>0$ such that $Z\leq \delta^3 C_2/a$.

Splitting the last term in (\ref{eqnDotVAux4}) into two sums depending on whether $|x_k|\leq a$ and substituting $Z$ for the sum with indices satisfying $|x_k|>a$, we can simplify  (\ref{eqnDotVAux4})  as
\begin{multline}\label{eqnVDotAux5}
\dot V  \leq    - F_{k_1}(x_{k_1})^T\cdot x_{x_1} + C_1|x_{k_1}| - \sum_{k\neq k_1} F_k(x_k)^T\cdot x_k \\
 -\frac{\delta}{a} \sum_{|x_k|\leq a} \left[y_k +F_k(x_k)- \delta W_k \right]^T\cdot  y_k + \delta^3 C_2,
\end{multline}

Since  the functions $p_k(|x|)$ are continuous and $p_k(|x|)|x|\to\infty $  as $|x|\to\infty$, there is a constant $C>0$
 such that $ p_k(|x_k|)  |x_k|^2 \geq -C$ for every $k=1,\ldots,n$ and every $x_k\in \mathbb R^d$. Hence, there exists
 $C_3>0$ such that $$\sum_{k\neq k_1} F_k(x_k)^T\cdot x_k = \sum_{k\neq k_1} p_k(|x_k|)|x_k|^2\geq -C_3$$ for any $k_1$ and  for any choice of the vectors $\{x_k\}$.   This implies that
\begin{multline}\label{eqnVDotAux6}
\dot V  \leq    - F_{k_1}(x_{k_1})^T\cdot x_{x_1} + C_1|x_{k_1}| +(C_3+ \delta^3 C_2)
 -\frac{\delta}{a} \sum_{|x_k|\leq a} \left[y_k +F_k(x_k)- \delta W_k \right]^T\cdot  y_k .
\end{multline}

\medskip
(IV) Now we are ready to select the parameters $a>1$, $0<\delta<\mbox{min}\{1, 1/||S||_{\infty}\}$, and $b>0$.  Choose $a>1$ large enough so that for every $k\in \{1,\ldots,n\}$, we have that
\begin{equation}\label{eqnChooseA}- p_k(z)z^2 +  C_1 z  + (C_3+ \delta^3 C_2) < -2 \mbox{ whenever }z>a.\end{equation}

 In particular, since we assumed that $ \delta $ is below $1/||S||_{\infty}$, we obtain  that $\delta|W_k|<1$. Set $D_k = F_k(x_k)-\delta W_k$. Given $y \in \mathbb R^d$, consider $[y+D_k]^T\cdot y$ as a function of $y$. It attains a global minimum at $y= -\frac{1}{2}D_k$. Therefore,  $[y+D_k]^T\cdot y \geq -\frac{1}{4}|D_k|^2$. Since $F_k(x_k)$ is bounded whenever $|x_k|\leq a$, $|W_k|\leq ||S||_{\infty}$ and $\delta<1$, we can find a constant $C_4>0$ (that depends on $a$) such that
 $$\left[y_k+ D_k\right]^T\cdot y_k \geq -C_4\mbox{ for all }k.$$

 Choose $0<\delta<1$ small enough so that \begin{equation}\label{eqnChooseDelta} \frac{\delta}{a} \sum_{|x_k|\leq a} \left[y_k +F_k(x_k)- \delta W_k \right]^T\cdot  y_k \geq -\frac{\delta n C_4}{a} >-1.\end{equation}

If $|x_k|>a $ for some index $k$, then $|x_{k_1}|>a$. Therefore, substituting  Inequalities (\ref{eqnChooseA}) and ({\ref{eqnChooseDelta}}) into (\ref{eqnVDotAux6}), we obtain that $\dot V \leq -2+1 = -1$.

\medskip
(III) The choice of $b$ is motivated by the need to address the case when $|x_k|\leq a$ for all $k$, but $|y_k|>b$ for some $k$.

 Set
$$W_a = \max_k\max_{|x|\leq a}\Big[- F_{k}(x)^T\cdot x + C_1|x| +(C_3+ \delta^3 C_2) \Big ].$$ It follows from (\ref{eqnVDotAux6}) that
\begin{equation}\label{eqnVDotAux7}
\dot V  \leq    W_a
 -\frac{\delta}{a} \sum_{|x_k|\leq a} \left[y_k +F_k(x_k)- \delta W_k \right]^T\cdot  y_k .
\end{equation}

 Fix an index $k_2$ such that $|y_{k_2}| = \max_k |y_k|$. Assume $|x_k|\leq a$ for every $k$. Then there exists  a constant  $K_a>0$ such that $|D_k|=|F_k(x_k)- \delta W_k|\leq K_a$ for any choice of vectors $\{x_1,\ldots,x_n\}$ provided that $|x_i|\leq a$.  Thus,  we can rewrite Inequality (\ref{eqnVDotAux7})  as
\begin{equation}\label{eqnVDotAux8}
\begin{split}\dot V & \leq W_a   -\frac{\delta}{a}  \left[y_{k_2}+ D_{k_2})\right]^T  \cdot y_{k_2} -\frac{\delta}{a}  \sum_{ k\neq k_2} \left[y_k+ D_k)\right]^T  \cdot y_k  \\
& \leq -\frac{\delta}{a}  \left[y_{k_2}+ D_{k_2}\right]^T  \cdot y_{k_2} + W_a + \frac{\delta}{4a}\sum_{ k\neq k_2} |D_k|^2   \\
& \leq -\frac{\delta}{a}\left [|y_{k_2}|^2 - |D_{k_2}|\cdot |y_{k_2}| \right] + W_a+\frac{\delta (n-1) K_a}{4a}
\\
& \leq -\frac{\delta}{a}\left [|y_{k_2}|^2 - K_a |y_{k_2}| \right] + W_a+\frac{\delta (n-1) K_a}{4a}.
\end{split}
\end{equation}

Choose $b>0$ large enough such that  $$\frac{\delta}{a}\left [z^2 - K_a z \right] -W_a-\frac{\delta (n-1) K_a}{4a} > 1\mbox{ whenever }z \geq b. $$

Thus, combining this inequality with (\ref{eqnVDotAux8}), we obtain that  $\dot V<-1$ whenever $|x_k|\leq a$ for all $k$, but $|y_k|>b$ for some $k$.

Thus, we have shown that the full weak derivative of the continuous function $V$  is less than $-1$ outside of the box $B$. Since the function $t\mapsto V(\gamma(t))$ is locally Lipschitz-continuous, it can be recovered from its weak derivative. Therefore, for any solution $\gamma = (x,y)$ that lies outside the box $B$ when $t\in [t_1,t_2]$, we have that
$$V(\gamma(t_2))-V(\gamma(t_1)) = \int_{t_1}^{t_2} \dot V(\gamma(t)) dt\leq  \int_{t_1}^{t_2} (-1) dt = - (t_2-t_1).$$

Notice also that $V\to\infty$ as $|x|+|y|\to\infty$. The result follows from the standard application of the Lyapunov method, see, for example, \cite[Theorem XVI]{LaSalle:Book:1961}.
$\hfill\square$
%

\begin{remark} In Theorem \ref{ThmMain} we can slightly weaken the assumption that the function $p_k$ is radial and instead assume that $p_k = p_k(x)$ depends directly on the velocity and there is a function $\phi_k: [0,\infty)\rightarrow \mathbb R$ such that $p_k(x)\geq \phi_k(|x|)$ and $\phi_k(z)z\to\infty$ as $z\to \infty$.
For example, the propelling capabilities of agents (say birds in $\mathbb R^3$) may be isotropic within the horizontal components, but different in the vertical direction.

\end{remark}
%

\bigskip
Consider the system
\begin{equation}\label{eqnParabolicPotentialModel} \ddot r_k = - p_k(|\dot r_k|)\dot r_k - \lambda (r_k-\bar r),\end{equation}
where $\lambda>0$,  $r_k\in \mathbb R^d$, and $\bar r= (r_1+\cdots + r_n)/n$, and $p_k(r)r\to\infty$ as $r\to\infty$. Setting $x_k = \dot r_k$, $y_k = \ddot r_k$, and applying Theorem \ref{ThmMain}, we immediately obtain the following result.

\begin{corollary}\label{CorollaryParabolicPotential} There are positive constants $C_1$, $C_2$, $C_3$, and $C_4$ such that for any solution $r(t)$ of System (\ref{eqnParabolicPotentialModel}) we have that $$|\dot r_k(t)|\leq C_1,\;|\ddot r_k(t)|\leq C_2,\;|r_k(t) - \bar r(t)|\leq C_3,\mbox{ and }|\bar r(t)| \leq C_4 t$$ for all $t$ large enough.
\end{corollary}

\begin{remark} A flock whose agents satisfy (\ref{eqnParabolicPotentialModel}) could have its center of mass move steadily away from the origin with   $|\bar r (t) |\geq C_5 t  $. Thus, the estimate given in  Corollary \ref{CorollaryParabolicPotential} is sharp.  To illustrate that, consider two agents in $\mathbf R^3$, propelled by a common function $p(|\dot r||)=1-|\dot r |^2.$ Let $a, b$ be positive such that $a^2+b^2=1.$ Then $r_1(t)= (at, b \cos ( \sqrt \lambda t), b \sin ( \sqrt \lambda t)^T, r_2(t)= (at, -b \cos ( \sqrt \lambda t), -b \sin ( \sqrt \lambda t)^T $ is a solution to (\ref{eqnParabolicPotentialModel}) that has $\bar r(t)= (at, 0, 0 )^T.$
\end{remark}

\begin{corollary}\label{CorollaryParabolicPotential2}
Consider the system
\begin{equation}\label{eqnParabolicPotentialModel2} \ddot r_k = - p_k(|\dot r_k|)\dot r_k - \sum _{m=1}^n a_{k,m} r_m,\end{equation}
where $A=\{a_{k,m}\}$ is symmetric and positive definite,  $r_k\in \mathbb R^d$, and $p_k(r)r\to\infty$ as $r\to\infty$.
There are positive constants $C_1$, $C_2$,and  $C_3$ such that for any solution $r(t)$ of System (\ref{eqnParabolicPotentialModel2}) we have that $$|\dot r_k(t)|\leq C_1,\;|\ddot r_k(t)|\leq C_2,\;|r_k(t) |\leq C_3$$ for all $t$ large enough.
\end{corollary}

Proof: The first two inequalities follow from setting $x_k = \dot r_k$, $y_k = \ddot r_k$, and applying Theorem \ref{ThmMain} for the collection of energies $  E_k=0 $ (recall that $Q=0$ when $A$ in invertible).
 The last inequality follows from
 $ r_k= \sum _{j=1}^n a_{kj}^{(-1)} \left [-y_j - p_j(|x_j|)x_j \right ] $ where $A^{-1}=\{a_{kj}^{(-1)}\}.$ $\hfill\square$

 \bigskip

We note that when the matrix $A$ is singular, the projection operator $Q:\mathbb R^n\rightarrow \textrm{Ker(A)}$ can be viewed as the generalized center of mass of the system.

\begin{corollary}\label{CorollaryParabolicPotential3}
Consider the system
(\ref{eqnParabolicPotentialModel2} )
where $A=\{a_{k,j}\}$ is symmetric, positive semidefinite, with projection matrix onto its kernel denoted $Q,$ and  $r_k\in \mathbb R^d$, and $p_k(r)r\to\infty$ as $r\to\infty$.
There are positive constants $C_1$, $C_2$,and  $C_3$ such that for any solution $r(t)$ of System (\ref{eqnParabolicPotentialModel2}) we have that $$|\dot r_k(t)|\leq C_1,\;|\ddot r_k(t)|\leq C_2,\;|r_k(t) - \sum _{j=1}^n q_{k,j} r_j(t)|\leq C_3$$ for all $t$ large enough.
\end{corollary}

Proof: The first two inequalities follow from setting $x_k = \dot r_k$, $y_k = \ddot r_k$, and applying Theorem \ref{ThmMain}.
Note that by multiplying the last equation of $(\ref{eqnParabolicPotentialModel2} )$ with $Q$ we conclude that all the constant energies $E_k$ associated with the system in $(x, y)$ are zero.
The last inequality uses  $ I-Q=(A+Q)^{-1} A$ to get
 $ r_k(t) - \sum _{j=1}^n q_{k,j} r_j(t)= \sum _{j=1}^n s_{kj}^{(-1)} \left [-y_j - p_j(|x_j|)x_j \right ] $ where $(A+Q)^{-1}=(s_{kj}^{(-1)}).$
$\hfill\square$

\bigskip 

\begin{figure}[h!]
\centering
\includegraphics[width=.7\textwidth]{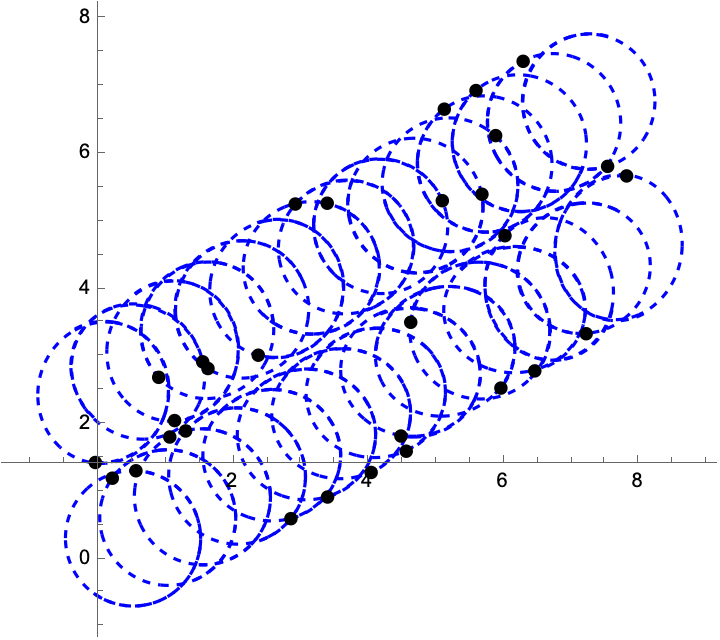}
\caption{\label{FigureRings} A simulation of a ring state for Equation  (\ref{eqnParabolicPotentialModel2}) with the communication topology as described in Example \ref{ExamplesGeomCenters}(2). The dots represent the agents and the circles represent their limit cycles.}
\end{figure}

 \begin{example}\label{ExamplesGeomCenters}
 
 (1)  If $A$ is a square matrix such that $\textrm{ker}(A) = \textrm{Span}\{(1,\ldots,1)^T\}$, then the projection operator $Q$ onto the kernel of $A$ is the averaging operator. Thus, in view of Corollary \ref{CorollaryParabolicPotential3}, the swarm will congregate around the center of mass of the system. Examples of such matrices include systems with communication topologies allowing the agents to communicate only with its left and right neighbors, say, symmetric nearest neighbor, and such that the row-sum of the weights is equal to zero. 
 
 (2)  
 Figure (\ref{FigureRings}) shows limit cycles for Equation (\ref{eqnParabolicPotentialModel2}) in $\mathbb R^2$ with $p_k(z) = (z^2-1)$ and $A = I-Q$, where $Q$ is the projection operator onto the subspace spanned by the vectors $(1,0,1,\ldots,0)^T$ and $(1/n,\ldots,n/n)^T$ for $n=30$ agents. The agents are given random uniformly-distributed positions in the square $[1,5]\times [1,5]$ and zero initial velocities. 

\end{example}

%
\section{Bounded Coupling}\label{SectionBoundedCoupling}
In this section, we establish the ultimate boundedness for systems with weak coupling.

\begin{theorem}\label{TheoremWeakCoupling}
Consider the system
\beq
 \ddot r_k = -p_k(|\dot r_k|) \dot r_k + c_k(r_1, r_2, \dots, r_n, \dot r_1, \dot r_2, \dots, \dot r_n),
 \label{BoundedProblem}
 \eeq
where the coupling functions $c_k:\mathbb R^{2dn}\rightarrow \mathbb R^{d}$ are bounded and locally Lipschits and the functions $p_k$ satisfy $p_k(r)\to \infty $ as $r \to \infty$. Then there exist constants $M_1$ and $ M_2$ such that for any solution $r=(r_1,\ldots,r_n)$ of (\ref{BoundedProblem}) we have that
 $|\dot r_k(t)|\leq M_1$ and  $|\ddot r_k(t)|\leq M_2$, for all large enough $t$.
\end{theorem}
{\it Proof.} To simplify the notation, we will simply write $c_k$ for the coupling functions. Choose $m>0$ such that $|c_k|=|c_k(r_1,\ldots,r_n,\dot r_1,\ldots \dot r_n)|\leq m$ for all $k$ and  $r_i$, $\dot r_i$. Note that
$$
|\ddot r_k| \leq |p_k(|\dot r_k|) |\dot r_k| + |c_k|\leq
|p_k(|\dot r_k|)||\dot r_k| +m.$$
Thus, the boundedness of the velocities automatically implies the boundedness of the accelerations.

For each $k$, consider the function from $[0, \infty )$ into $\R$ given by $s \to s[p_k(s)-m].$ Note that $s[p_k(s)-m]\to\infty$ as $s\to \infty$. Therefore, there is  a global minimum $q_k$ (necessarily $q_k\leq 0$). That is, $s[p_k(s)-m]\geq q_k$ for all $s\geq 0$. Set $Q=-\sum_{k} q_k. $ Then for any index $k_1$ and any $s\geq 0$ we get
\beq
\sum _{ k \neq k_1} [p_k(s)s^2-ms]\geq -Q
\label{lowboundB}
\eeq
Let $M_1$ be large enough such that for any $k=1,\ldots,n$, we have  that $ s[p_k(s)-m] \geq Q+1$ whenever $s\geq M_1$.

Consider the kinetic energy of the system
$E =\frac{1}{2} \sum _{k} \dot r_k^T \cdot \dot r_k$. Assume that $|\dot r_{k_1}|\geq M_1$ for some index $k_1$. Differentiating the function $E$ along  trajectories of the system, we obtain that
\begin{equation*}\begin{split}\dot E &= \sum_k \big [ -p_k(|\dot r_k|)|\dot r_k|^2+ \dot r_k^T \cdot c_k \big ] \\
&\leq \sum _k \big[ -p_k(|\dot r_k|)|\dot r_k|^2+m|\dot r_k|  \big] \\
& = - \big[p_{k_1}(|\dot r_k|)|\dot r_{k_1}|-m\big]|\dot r_k|+\sum _{k\neq k_1} \big[ -p_k(|\dot r_k|)|\dot r_k|+m|\dot r_k|  \big] \\
& \leq -(Q+1)- \sum _{k\neq k_1} \big[ p_k(|\dot r_k|)|\dot r_k|-m|\dot r_k|  \big].
\end{split}
\end{equation*}
Combining this inequality with (\ref{lowboundB}), we obtain that $\dot E\leq -1$ whenever $|\dot r_k|\geq M_1$ for at least one index $k$. Thus, using the function $E$ as a  Lyapunov function for the system, we can conclude that the velocities are ultimately bounded \cite[Theorem XVI]{LaSalle:Book:1961}.

$\hfill\square$

\begin{remark} A flock whose agents satisfy the assumptions of Theorem \ref{TheoremWeakCoupling} could disperse, in the sense that $|r_1(t)-\bar r (t)|$ may go to infinity as $t\to \infty.$
Consider for example the two-agent flock in $\mathbf R$ governed by
$$ \ddot x_1= (1-\dot x _1^2)\dot x _1 + f(x_1, x_2, \dot x_1, \dot x_2), \; \; \ddot x_2= (1-\dot x _2^2)\dot x _2 - f(x_1, x_2, \dot x_1, \dot x_2)$$
where the $C^1$ function $f$ satisfies: there exists $M$ such that  $|f(x_1, x_2, v_1, v_2) |< m/2 $ for all $(x_1, x_2, v_1, v_2)$ with $\mbox{max} \{|x_1|, |x_2| \}\geq M $ where $-m$ denotes the minimum of cubic polynomial $v^3-v$ for $v \in[0, \infty) $, $m \approx 0.385$ (attained when $v=1/\sqrt 3$).

Consider the preimage of $[-m/2, m/2]$ under the polynomial $v^3-v;$ denote the component containing 1 by $[v_a, v_b],  $ thus $[v_a, v_b]\approx [0.885, 1.085].  $

Two agents starting with initial conditions opposite from each other remain opposite at all times, with $ x(t)=x_1(t) =-x_2(t)$ satisfying $$ \ddot x= (1-\dot x ^2)\dot x -f(x, -x, \dot x, -\dot x).$$
Start with $x(0)>M$ and $\dot r(0) $ between  $v_a$ and   $v_b$. One can show that for all $t>0$ we have $x(t)>M$ and $v_a< \dot r(t)< v_b, $ therefore $x(t) > M+v_at.$ (Assume, by contradiction, that there is a first time $t_1>0$ when $x(t_1)\leq M$ or $\dot x(t_1)\leq v_a$ or $\dot x(t_1)\geq v_b.$ It must be that $\dot x(t_1)=v_a$ or $\dot x(t_1)=v_b$ since while the velocity is in $[v_a, v_b]$ the position is increasing. Note that when $x>M, \dot x=v_a$, $\ddot x =m/2+f(x, -x, v_a, -v_a)>0$ so $\dot x $ is strictly increasing, meaning that $t_1$ could not have been the first time $\dot x$ dropped below $v_a.$ Similarly, when $x>M, \dot x=v_b$ we have $\ddot x <0 $ so the velocity is decreasing, making it impossible for $t_1$ to be the first time $\dot x \geq v_b.$ )

\end{remark}

\printbibliography

\end{document}